\documentstyle[a4,11pt]{article}
\begin{document}
\begin{titlepage}
\vskip 2cm
\begin{flushright}
Preprint CNLP-1994-11
\end{flushright}
\vskip 2cm
\begin{center}
{\bf
THE L-EQUIVALENT COUNTERPART OF THE M-III EQUATION}
\footnote{Preprint
CNLP-1994-11. Alma-Ata. 1994 }
\end{center}
\vskip 2cm
\begin{center}
R. Myrzakulov
\footnote{E-mail: cnlpmyra@satsun.sci.kz}
and A. K. Danlybaeva
\footnote{E-mail: cnlpakd@satsun.sci.kz}
\end{center}
\vskip 1cm

\begin{center}
 Centre for Nonlinear Problems, PO Box 30, 480035, Alma-Ata-35, Kazakhstan
\end{center}

\begin{abstract}
The connection between differential geometry of curves and the integrable (2+1)
-dimensional spin system (the M-III equation) is established. Using the
proposed geometrical formalism, the L-equivalent counterpart of the M-III
equation is found.
\end{abstract}

\end{titlepage}

\setcounter{page}{1}
\newpage

\tableofcontents
\section{Introduction}

Consider the Myrzakulov III ( M-III ) equation [1]
$$
{\bf S}_{t}  =  ({\bf S}\wedge {\bf S}_{y}+u{\bf S})_{x} +
2l(cl+d){\bf S}_{y} -4cv{\bf S}_{x}  + {\bf S}\wedge {\bf W} \eqno (1a)
$$
$$
u_{x} = - {\bf S}\cdot ({\bf S}_{x}\wedge {\bf S}_{y})
   \eqno (1b)
$$
$$
v_{x} = \frac{1}{4(2cl+d)^{2}}({\bf S}^{2}_{x})_{y}   \eqno(1c)
$$
$$
W_{x} = J {\bf S}_{y}   \eqno(1d)
$$
which was introduced in [1]. Here $ {\bf S}=(S_1,S_2,S_3) $ is the spin vector,
$ {\bf S}^2=\beta=\pm1, $ $u$ and $v$ are scalar functions, $c, d, l$ - constants.
The M-III equation (1) contains several interesting integrable particular
cases:

i) $c=0,d=1$,   yields the isotropic Myrzakulov I ( M-I ) equation
$$
{\bf S}_{t}  =  ({\bf S}\wedge {\bf S}_{y}+u{\bf S})_{x}   \eqno (2a)
$$
$$
u_{x} = - {\bf S}\cdot ({\bf S}_{x}\wedge {\bf S}_{y}).   \eqno (2b)
$$

ii) $d=0$ yields  the Myrzakulov II ( M-II ) equation
$$
{\bf S}_{t}  =  ({\bf S}\wedge {\bf S}_{y}+u{\bf S})_{x} +
2cl^{2}{\bf S}_{y} -4cv{\bf S}_{x}  \eqno (3a)
$$
$$
u_{x} = - {\bf S}\cdot ({\bf S}_{x}\wedge {\bf S}_{y})   \eqno (3b)
$$
$$
v_{x} = \frac{1}{16c^{2}l^{2}}({\bf S}^{2}_{x})_{y}.   \eqno(3c)
$$

iii) $c=0, d=1, J=diag(0,0,\triangle)$,  yields the M-I equation with one-ion
anisotropy
$$
{\bf S}_{t}  =  ({\bf S}\wedge {\bf S}_{y}+u{\bf S})_{x} +
2l{\bf S}_{y}   + {\bf S}\wedge {\bf W} \eqno (4a)
$$
$$
u_{x} = - {\bf S}\cdot ({\bf S}_{x}\wedge {\bf S}_{y})   \eqno (4b)
$$
$$
W_{x} = J {\bf S}_{y}.   \eqno(4c)
$$

iv) The isotropic M-III equation as $ J=diag(0,0,0)$
$$
{\bf S}_{t}  =  ({\bf S}\wedge {\bf S}_{y}+u{\bf S})_{x} +
2l(cl+d){\bf S}_{y} -4cv{\bf S}_{x} \eqno (5a)
$$
$$
u_{x} = - {\bf S}\cdot ({\bf S}_{x}\wedge {\bf S}_{y})   \eqno (5b)
$$
$$
v_{x} = \frac{1}{4(2cl+d)^{2}}({\bf S}^{2}_{x})_{y}   \eqno(5c)
$$
and so on [1].    All of these equations are integrable. For instance,
the Lax representation of the isotropic M-III equation (5) has  the   form [1]
$$
\phi_{x} = U' \phi   \eqno(6a)
$$
$$
\phi_{t} = -2(c\lambda^2+d\lambda)\Phi_y+ V' \phi   \eqno(6b)
$$
with
$$
U' = [ic(\lambda^2-l^2)+id(\lambda-l)]S+\frac{c(\lambda-l)}{2c\lambda+d}SS_x,
\eqno(7a)
$$
$$
V' =[2c(\lambda^2-l^2)+2d(\lambda-l)]B+\lambda^2 F_2+\lambda F_1+F_0,
         \eqno(7b)
$$
where
$$
F_2=-4ic^2VS,
$$
$$
F_1=-4icdVS-\frac{4c^2 V}{2cl+d}VSS_x-\frac{ic}{2cl+d}S{(SS_x)_y-[SS_x,B]},
$$
$$
F_0=-lF_1-l^2F_2, B=\frac{1}{4}([S,S_y]+2iuS), S=\vec S \cdot \vec \sigma.
$$

In [1] we proposed a new class integrable and nonintegrable spin systems. And
we also suggested the geometrical formalism to establish the connection
between differential geometry of curves and surfaces and nonlinear evolution
equations (NEE), including soliton equations (the A-, B-, C-, D-approaches)
(see, also, the refs. [4, 5]).
In this paper, using the D-approach we will establish the connection
between the differential geometry of curves and the isotropic  M-III
equation (5). Also we find the L-equivalent (Lakshmanan equivalent [1, 2])
counterpart of this equation.

\section{On the (2+1)-dimensional curve  and soliton equations.}
\subsection{The M-LXVII equation}

Using D-apprach, in [1]  we have  established
the connection between differential geometry of  curves and well known
soliton equations in 2+1 dimensions. It is remarkable that this approach simultaneously permit determine the
L-equivalent counterpart of the under consideration spin systems. Here we
will consider a (2+1)-dimensional curves which are
given by the M-LXVII  equation. The M-LXVII  equation reads as [1]
$$
\left ( \begin{array}{c}
{\bf e}_{1} \\
{\bf e}_{2} \\
{\bf e}_{3}
\end{array} \right)_{x} =
\left ( \begin{array}{ccc}
0             & k     & -\gamma \\
-\beta k      & 0     & \tau  \\
\beta\gamma   & -\tau & 0
\end{array} \right)
\left ( \begin{array}{ccc}
{\bf e}_{1} \\
{\bf e}_{2} \\
{\bf e}_{3}
\end{array} \right) \eqno(8a)
$$
$$
\left ( \begin{array}{ccc}
{\bf e}_{1} \\
{\bf e}_{2} \\
{\bf e}_{3}
\end{array} \right)_{t}= C
\left ( \begin{array}{ccc}
{\bf e}_{1} \\
{\bf e}_{2} \\
{\bf e}_{3}
\end{array} \right)_{y} +  D
\left ( \begin{array}{ccc}
{\bf e}_{1} \\
{\bf e}_{2} \\
{\bf e}_{3}
\end{array} \right) \eqno(8b)
$$
where we assume that  $k=k(\lambda, x, y,t), \quad \tau =
\tau(\lambda, x, y, t), \quad \gamma = \gamma (\lambda, x,y,t), \quad
C = C(\lambda, x, y, t),$ \\
$  D = D(\lambda, x, y, t)$ and $\lambda$ is
some complex parameter, $C, D$ are some matrices. And also let us we
suppose that
$$
k = \sum^{n}_{j=1} k_{j}\lambda^{j}, \quad
\gamma = \sum^{n}_{j=1} \gamma_{j}\lambda^{j}, \quad
\tau = \sum^{n}_{j=1} \tau_{j}\lambda^{j}  \eqno(9a)
$$
$$
C = \sum^{m}_{j=1} C_{j}\lambda^{j}, \quad
D = \sum^{m}_{j=1} D_{j}\lambda^{j},\quad e_1^2=\beta=\pm1,\quad e_2^2=e_3^2=1.
\eqno(9b)
$$
Here $k_{j}=k_{j}(x,y,t), \quad \tau_{j} =
\tau_{j}(x, y, t), \quad \gamma_{j} = \gamma_{j}(x,y,t), \quad
C_{j} = C_{j}(x, y, t), \quad D_{j} = D_{j}(x, y, t)$.

\subsection{The M-LXVII equation associated with the M-III equation}

In this subsection we require that the vector ${\bf e}_{1} \equiv {\bf S}$ satisfies
the M-III equation (5). Then we must put
$$
k_{0} =  -d(q+p), \quad k_{1} = -2c(q+p), \quad k_{j} = 0, \quad j\ge 2
\eqno(10a)
$$
$$
\gamma_{0} =  -d(q-p), \quad \gamma_{1} = -2c(q-p), \quad
\gamma_{j} = 0, \quad  j\ge 2
\eqno(10b)
$$
$$
\tau_{0} = 0, \quad \tau_{1} = -2d, \quad \tau_{2} = -2c,
\quad \tau_{j} = 0,   \quad  j\ge 3
\eqno(10c)
$$
$$
C_{0} = 0, \quad C_{1} = 2d, \quad C_{2} = 2c, \quad C_{j} = 0  \quad  j\ge 3
\eqno(10d)
$$
$$
D = \lambda^{2} D_{2} + \lambda D_{1} + D_{0},
\quad D_{j} = 0  \quad  j\ge 3.
\eqno(10e)
$$
Here
$$
D=\lambda^2 D_2+\lambda D_1+D_0 =
\left ( \begin{array}{ccc}
0       & \omega_{3}  & -\omega_{2} \\
-\beta\omega_{3} & 0      & \omega_{1} \\
\beta\omega_{2}  & -\omega_{1} & 0
\end{array} \right)  \eqno(11)
$$
$q, p$ are some function. So  in our case the M-LXVII equation takes the form [1]
$$
\left ( \begin{array}{c}
{\bf e}_{1} \\
{\bf e}_{2} \\
{\bf e}_{3}
\end{array} \right)_{x} =
\left ( \begin{array}{ccc}
0             & -(2c\lambda + d)(q+p)     & i(2c\lambda + d)(q-p) \\
\beta(2c\lambda+d)(q+p)      & 0     & -2(c\lambda^{2}+d\lambda)  \\
-\beta i(2c\lambda+d)(q-p)   & 2(c\lambda^{2}+d\lambda) & 0
\end{array} \right)
\left ( \begin{array}{ccc}
{\bf e}_{1} \\
{\bf e}_{2} \\
{\bf e}_{3}
\end{array} \right) \eqno(12a)
$$
$$
\left ( \begin{array}{ccc}
{\bf e}_{1} \\
{\bf e}_{2} \\
{\bf e}_{3}
\end{array} \right)_{t}= 2(c\lambda^{2} + d \lambda)
\left ( \begin{array}{ccc}
{\bf e}_{1} \\
{\bf e}_{2} \\
{\bf e}_{3}
\end{array} \right)_{y} +  (\lambda^{2}D_{2}+\lambda D_{1} + D_{0})
\left ( \begin{array}{ccc}
{\bf e}_{1} \\
{\bf e}_{2} \\
{\bf e}_{3}
\end{array} \right) \eqno(12b)
$$
where the explicit forms of $ D_j $ are given in [1].
In terms of matrix this equation we can write in the form
$$
\hat e_{1x}  = -(2c\lambda + d)(q+p)\hat e_{2}+i(2c\lambda+d)(q-p)\hat e_{3}  \eqno(13a)
$$
$$
\hat e_{2x}  =\beta (2c\lambda + d)(q+p)\hat  e_{1}-2(c\lambda^{2}+d\lambda)\hat e_{3}
\eqno(13b)
$$
$$
\hat e_{3x}  =-\beta
 i(2c\lambda + d)(q-p)\hat  e_{1}+2(c\lambda^{2}+d\lambda)\hat e_{2}
\eqno(13c)
$$
$$
\hat e_{1t} =  2(c\lambda^{2}+d\lambda) \hat e_{1y} +\omega_{3} \hat e_{2}-
\omega_{2}\hat e_{3}  \eqno(13d)
$$
$$
\hat e_{2t} =  2(c\lambda^{2}+d\lambda) \hat e_{2y}
  - \beta \omega_{3} \hat e_{1} +\omega_{1}\hat e_{3}  \eqno(13e)
$$
$$
\hat e_{3t} =  2(c\lambda^{2}+d\lambda) \hat e_{3y}
+\beta \omega_{2} \hat e_{1}-\omega_{1}\hat e_{2}  \eqno(13f)
$$
where
$$
\hat  e_{1} = g^{-1}\sigma_{3}g, \quad \hat e_{2}=g^{-1}\sigma_{2}g,
\quad \hat e_{3} = g^{-1}\sigma_{1}g.  \eqno(14)
$$
Here  $\sigma_{j}$ are Pauli matrices
$$
\sigma_{1} =
\left ( \begin{array}{cc}
0       &  1 \\
1       &  0
\end{array} \right) , \quad
\sigma_{2} =
\left ( \begin{array}{cc}
0       &  -i \\
i       &  0
\end{array} \right) , \quad
\sigma_{3} =
\left ( \begin{array}{cc}
1       &  0 \\
0       &  -1
\end{array} \right).   \eqno(15)
$$
So we have
$$
\sigma_{1}\sigma_{2} = i\sigma_{3} = -\sigma_{2}\sigma_{1}, \quad
\sigma_{1}\sigma_{3} = -i\sigma_{2} = -\sigma_{3}\sigma_{1}, \quad
\sigma_{3}\sigma_{2} = -i\sigma_{1} = -\sigma_{2}\sigma_{3}
\eqno(16a)
$$
and
$$
\sigma_{j}^{2} =I = diag(1,1). \eqno(16b)
$$
Equations (13) we can rewrite in the form
$$
[\sigma_{3}, U] = -(2c\lambda + d)(q+p)\sigma_{2}+i(2c\lambda+d)(q-p)\sigma_{1}  \eqno(17a)
$$
$$
[\sigma_{2}, U] =
\beta (2c\lambda + d)(q+p)\sigma_{3}-2(c\lambda^{2}+d\lambda)\sigma_{1}
 \eqno(17b)
$$
$$
[\sigma_{1}, U] =
 -\beta i(2c\lambda + d)(q-p)\sigma_{3}+2(c\lambda^{2}+d\lambda)\sigma_{2}
 \eqno(17c)
$$
$$
[\sigma_{3}, V] =  \omega_{3} \sigma_{2}-
\omega_{2}\sigma_{1}  \eqno(17d)
$$
$$
[\sigma_{2}, V] =  -\beta \omega_{3} \sigma_{3}
+\omega_{1}\sigma_{1}  \eqno(17e)
$$
$$
[\sigma_{1}, V] = \beta \omega_{2} \sigma_{3}
-\omega_{1}\sigma_{2}.  \eqno(17f)
$$
Here
$$
U = g_{x}g^{-1}, \quad  V = g_{t}g^{-1} - 2(c\lambda^{2}+d\lambda)g_{x}g^{-1}.
\eqno(18)
$$
Hence we get
$$
U = i[(c\lambda^{2}+d\lambda)\sigma_{3}+(2c\lambda+d)Q]    \eqno(19a)
$$
$$
V = \lambda^{2}B_{2}+\lambda B_{1}+B_{0}  \eqno(19b)
$$
with
$$
B_{2} = -4ic^{2}v\sigma_{3}, \quad
B_{1} = -4icdv\sigma_{3}+2cQ_{y}\sigma_{3}-8ic^{2}vQ, \quad
B_{0} = \frac{d}{2c}B_{1}-\frac{d^{2}}{4c^{2}}B_{2}  \eqno(20)
$$
and
$$
Q =
\left ( \begin{array}{cc}
0       &  q \\
p       &  0
\end{array} \right), \quad p=\beta \bar q.   \eqno(21)
$$
Thus the matrix-function  $g$ satisfies the equations
$$
g_{x} = Ug  \eqno(22a)
$$
$$
g_{t} = 2(c\lambda^{2}+d\lambda)g_{y}+Vg \eqno(22b)
$$
The compatibility condition of  these equations  gives the following M-III$_{q}$
equation
$$
iq_{t}=  q_{xy} -4ic(vq)_{x}+2d^{2}vq = 0    \eqno(23a)
$$
$$
-ip_{t}=  p_{xy} +4ic(vp)_{x}+2d^{2}vp = 0    \eqno(23b)
$$
$$
v_{x}= (pq)_{y}    \eqno(23c)
$$

So we have identified the curve,
given by the M-LXVII equation (12) with
the M-III equation (5). On the other hand, the compatibilty
condition of equations (22) is equivalent to the equation (23).
So that we have also established the connection between the  curve
(the M-LXVII  equation) and the equation (23).  It means that
the M-III equation (5) and the equation (23) are L-equivalent (Lakshmanan equivalent)
to each other.

\subsection{Particular cases}
\subsubsection{The M-LXVII equation and the M-I equation}
For the M-I equation $(c=0, d=1)$ the associated  M-LXVII equation has the form
$$
\left ( \begin{array}{c}
{\bf e}_{1} \\
{\bf e}_{2} \\
{\bf e}_{3}
\end{array} \right)_{x} =
\left ( \begin{array}{ccc}
0             & -(q+p)     & i(q-p) \\
\beta (q+p)      & 0     & -2\lambda  \\
-\beta i(q-p)   & 2\lambda & 0
\end{array} \right)
\left ( \begin{array}{ccc}
{\bf e}_{1} \\
{\bf e}_{2} \\
{\bf e}_{3}
\end{array} \right) \eqno(24a)
$$
$$
\left ( \begin{array}{ccc}
{\bf e}_{1} \\
{\bf e}_{2} \\
{\bf e}_{3}
\end{array} \right)_{t}= 2 \lambda
\left ( \begin{array}{ccc}
{\bf e}_{1} \\
{\bf e}_{2} \\
{\bf e}_{3}
\end{array} \right)_{y} +   D_{0}
\left ( \begin{array}{ccc}
{\bf e}_{1} \\
{\bf e}_{2} \\
{\bf e}_{3}
\end{array} \right) \eqno(24b)
$$
where  the explicit form  of $ D_0 $ given in [1].
 After some algebra we obtain the following L-equivalent of
the isotropic M-I equation (2)
$$
iq_{t}=  q_{xy} +2vq     \eqno(26a)
$$
$$
-ip_{t}=  p_{xy} +2vp     \eqno(26b)
$$
$$
v_{x}= (pq)_{y}    \eqno(26c)
$$
It is the Zakharov equation (ZE) [6].

\subsubsection{The M-LXVII equation and the M-II equation}

Now we put $d=0$. Then we get the following version of the M-LXVII
equation
$$
\left ( \begin{array}{c}
{\bf e}_{1} \\
{\bf e}_{2} \\
{\bf e}_{3}
\end{array} \right)_{x} =
\left ( \begin{array}{ccc}
0             & -2c\lambda (q+p)     & 2ci\lambda (q-p) \\
2\beta c\lambda (q+p)      & 0     & -2c\lambda^{2}  \\
-2\beta ci\lambda (q-p)   & 2c\lambda^{2} & 0
\end{array} \right)
\left ( \begin{array}{ccc}
{\bf e}_{1} \\
{\bf e}_{2} \\
{\bf e}_{3}
\end{array} \right) \eqno(27a)
$$
$$
\left ( \begin{array}{ccc}
{\bf e}_{1} \\
{\bf e}_{2} \\
{\bf e}_{3}
\end{array} \right)_{t}= 2c\lambda^{2}
\left ( \begin{array}{ccc}
{\bf e}_{1} \\
{\bf e}_{2} \\
{\bf e}_{3}
\end{array} \right)_{y} +  (\lambda^{2}D_{2}+\lambda D_{1} + D_{0})
\left ( \begin{array}{ccc}
{\bf e}_{1} \\
{\bf e}_{2} \\
{\bf e}_{3}
\end{array} \right) \eqno(27b)
$$
where the explicit forms of $ D_j $ are given in [1].
This M-LXVII equation is associated with the M-II equation (3). Proceeding as above
we get the following L-equivalent of the M-II equation
$$
iq_{t}=  q_{xy} -4ic(vq)_{x}    \eqno(28a)
$$
$$
-ip_{t}=  p_{xy} +4ic(vp)_{x}     \eqno(28b)
$$
$$
v_{x}= (pq)_{y}    \eqno(28b)
$$
which is the Strachan equation [7].

\section{The M-LXI equation and    the M-III equation}

\subsection{The M-LXI equation}

Now,  we consider the (2+1)-dimensional curve which is given by the
M-LXI equation
$$
\left ( \begin{array}{c}
{\bf e}_{1} \\
{\bf e}_{2} \\
{\bf e}_{3}
\end{array} \right)_{x}= A
\left ( \begin{array}{ccc}
{\bf e}_{1} \\
{\bf e}_{2} \\
{\bf e}_{3}
\end{array} \right) \eqno(29a)
$$
$$
\left ( \begin{array}{ccc}
{\bf e}_{1} \\
{\bf e}_{2} \\
{\bf e}_{3}
\end{array} \right)_{y}= B
\left ( \begin{array}{ccc}
{\bf e}_{1} \\
{\bf e}_{2} \\
{\bf e}_{3}
\end{array} \right)
\eqno(29b)
$$
$$
\left ( \begin{array}{ccc}
{\bf e}_{1} \\
{\bf e}_{2} \\
{\bf e}_{3}
\end{array} \right)_{t}= C
\left ( \begin{array}{ccc}
{\bf e}_{1} \\
{\bf e}_{2} \\
{\bf e}_{3}
\end{array} \right) \eqno(29c)
$$
where
$$
A =
\left ( \begin{array}{ccc}
0             & k     &  0 \\
-\beta k      & 0     & \tau  \\
0             & -\tau & 0
\end{array} \right) \eqno(30a)
$$
$$
B =
\left ( \begin{array}{ccc}
0            & m_{3}  & -m_{2} \\
-\beta m_{3} & 0      & m_{1} \\
\beta m_{2}  & -m_{1} & 0
\end{array} \right)  \eqno(30b)
$$
$$
C =
\left ( \begin{array}{ccc}
0       & \omega_{3}  & -\omega_{2} \\
-\beta\omega_{3} & 0      & \omega_{1} \\
\beta\omega_{2}  & -\omega_{1} & 0
\end{array} \right)  \eqno(30c)
$$

\subsection{The M-LXII equation}

From (29), we obtain the following M-LXII equation [1]
$$
A_y - B_x + [A, B] = 0 \eqno (31a)
$$
$$
A_t - C_x + [A, C] = 0 \eqno (31b)
$$
$$
B_t - C_y + [B, C] = 0 \eqno (31c)
$$
Equation (31a) gives
$$
k_{y} - m_{3x}  - \tau m_{2} = 0 \eqno(32a)
$$
$$
 - m_{2x} + \tau m_{3} - km_{1} =0   \eqno(32b)
$$
$$
\tau_{y} - m_{1x} + \beta km_{2}  =0. \eqno(32c)
$$

The M-LXII  equation (31), we can rewrite in form
$$
k_{y} - m_{3x} = \beta {\bf e}_{3}\cdot ({\bf e}_{3x}\wedge {\bf e}_{3y}) \eqno(33a)
$$
$$
 - m_{2x} = \beta {\bf e}_{2}\cdot ({\bf e}_{2x}\wedge {\bf e}_{2y}) \eqno(33b)
$$
$$
\tau_{y} - m_{1x}  = {\bf e}_{1}\cdot ({\bf e}_{1x}\wedge {\bf e}_{1y}) \eqno(33c)
$$

Also from (31) we get
$$
k_{t} - \omega_{3x} = \tau \omega_{2} = \beta  {\bf e}_{3}\cdot ({\bf e}_{3x}\wedge {\bf e}_{3t})
 \eqno (34a)
$$
$$
 - \omega_{2x} = - \tau \omega_{3} + k \omega_{1} =  \beta {\bf e}_{2}\cdot ({\bf e}_{2x}\wedge {\bf e}_{2t})
 \eqno (34b)
$$
$$
\tau_{t} - \omega_{1x} =- \beta k \omega_{2}  = {\bf e}_{1}\cdot ({\bf e}_{1x}\wedge {\bf e}_{1t})
 \eqno (34c)
$$
and
$$
m_{1t}-  \omega_{1y}= - \beta (m_{3} \omega_{2} - m_{2} \omega_{3}) =
 {\bf e}_{1}\cdot ({\bf e}_{1y}\wedge {\bf e}_{1t}) \eqno (35a)
$$
$$
m_{2t} - \omega_{2y} =- m_{1} \omega_{3} + m_{3} \omega_{1} = \beta
{\bf e}_{2}\cdot ({\bf e}_{2y}\wedge {\bf e}_{2t})
 \eqno (35b)
$$
$$
m_{3t} - \omega_{3y} =- m_{2} \omega_{1} + m_{1} \omega_{2} =  \beta
 {\bf e}_{3}\cdot ({\bf e}_{3y}\wedge {\bf e}_{3t})
 \eqno (35c)
$$

\subsection{On the topological invariants}

It is interesting to note that the M-LXII equations allows the following
integrals of motion [1]
$$
K_{1} = \int \int \kappa m_2 dxdy, \quad
K_{2} = \int \int \tau m_2 dxdy, \quad
K_{3} = \int \int (\tau m_3 -km_{1})dxdy \eqno(36a)
$$
or
$$
K_{1} = \int \int {\bf e}_1({\bf e}_{1x} \wedge {\bf e}_{1y})dxdy \eqno(37b)
$$
$$
K_{2} = \int \int {\bf e}_2({\bf e}_{2x} \wedge {\bf e}_{2y})dxdy \eqno(37b)
$$
$$
K_{3} = \int \int {\bf e}_3({\bf e}_{3x} \wedge {\bf e}_{3y})dxdy. \eqno(37c)
$$

So we have the following three topological invariants
$$
Q_{1} =\frac{1}{4\pi} \int \int {\bf e}_{1}\cdot ({\bf e}_{1x}\wedge {\bf e}_{1y})dxdy
\eqno(38a)
$$
$$
Q_{2} =\frac{1}{4\pi} \int \int {\bf e}_{2}\cdot ({\bf e}_{2x}\wedge {\bf e}_{2y})dxdy
\eqno(38b)
$$
$$
Q_{3} =\frac{1}{4\pi} \int \int {\bf e}_{3}\cdot ({\bf e}_{3x}\wedge {\bf e}_{3y})dxdy
\eqno(38c)
$$

\subsection{The M-III equation}

Now we wish show how connected  the M-LXI (29), M-LXII (31) and M-III
(5) equations. Let us we identify $ {\bf  S}$ with ${\bf e}_{1}$, i.e.
$$
{\bf e}_{1}  \equiv {\bf S}.   \eqno(39)
$$
Then we have
$$
m_{1}=u+\partial_{x}^{-1}\tau_{y} \eqno(40a)
$$
$$
m_{2}=\frac{1}{k} u_x  \eqno (40b)
$$
$$
m_{3} =\partial_x^{-1}(k_y-\tau m_2)=\partial^{-1}_x-\frac{\tau u_x}{k}
 \eqno(40c)
$$
and
$$
\omega_{1} = \frac{1}{k}[-\omega_{2x}+\tau\omega_{3}]
                 \eqno (41a)
$$
$$
\omega_{2}= -m_{3x}-\tau m_2+2l(cl+d)m_2
\eqno (41b)
$$
$$
\omega_{3}= m_{2x}-\tau m_3-uk+2l(cl+d)m_3-4cvk.
\eqno (41c)
$$

\subsection{The L-equivalent of the M-III equation}
We now  introduce the function $q$ according to the following expression
$$
q =\frac{k}{2(2cl+d)}exp i\{2l(cl+d)x-\partial^{-1}_x \tau\}.  \eqno(42)
$$

It is not difficult to verify that the
function $q$ satisfy the following nonlinear Schrodinger - type equation
$$
iq_{t}=  q_{xy} -4ic(vq)_{x}+2d^{2}vq = 0    \eqno(43a)
$$
$$
-ip_{t}=  p_{xy} +4ic(vp)_{x}+2d^{2}vp = 0    \eqno(43b)
$$
$$
v_{x}= (pq)_{y}.    \eqno(43c)
$$
It coincide wih the $ M-III_q $ equation (23). So we have proved that the equations
(23) and (5) are L-equivalent to
each other. As well known that these equations are G-equivalent to
each other [3]. Equation (43)  contains two reductions:
the Zakharov equation  as $c = 0$ [6] and
the Strachan  equation as $d = 0$ [7].

\section{The mM-LXI equation and the M-III equation}

\subsection{The mM-LXI equation}

One of the significant model of the (2+1)-dimensional curves is the
modified M-LXI (mM-LXI) equation. It reads  as [1]
$$
\left ( \begin{array}{c}
{\bf e}_{1} \\
{\bf e}_{2} \\
{\bf e}_{3}
\end{array} \right) = A_{m}
\left ( \begin{array}{ccc}
{\bf e}_{1} \\
{\bf e}_{2} \\
{\bf e}_{3}
\end{array} \right) \eqno(44a)
$$
$$
\left ( \begin{array}{ccc}
{\bf e}_{1} \\
{\bf e}_{2} \\
{\bf e}_{3}
\end{array} \right)_{y}= B_{m}
\left ( \begin{array}{ccc}
{\bf e}_{1} \\
{\bf e}_{2} \\
{\bf e}_{3}
\end{array} \right)
\eqno(44b)
$$
$$
\left ( \begin{array}{ccc}
{\bf e}_{1} \\
{\bf e}_{2} \\
{\bf e}_{3}
\end{array} \right)_{t}= C_{m}
\left ( \begin{array}{ccc}
{\bf e}_{1} \\
{\bf e}_{2} \\
{\bf e}_{3}
\end{array} \right) \eqno(44c)
$$
where
$$
A_{m} =
\left ( \begin{array}{ccc}
0             & k     & -\sigma \\
-\beta k      & 0     & \tau  \\
\beta\sigma   & -\tau & 0
\end{array} \right) \eqno(45a)
$$
$$
B_{m} =
\left ( \begin{array}{ccc}
0            & m_{3}  & -m_{2} \\
-\beta m_{3} & 0      & m_{1} \\
\beta m_{2}  & -m_{1} & 0
\end{array} \right)  \eqno(45b)
$$
$$
C_{m} =
\left ( \begin{array}{ccc}
0       & \omega_{3}  & -\omega_{2} \\
-\beta\omega_{3} & 0      & \omega_{1} \\
\beta\omega_{2}  & -\omega_{1} & 0
\end{array} \right).  \eqno(45c)
$$

\subsection{The mM-LXII equation}

The compatibility condition of the mM-LXI equations (44) gives the
modified M-LXII (mM-LXII) equation [1]
$$
A_{my} - B_{mx} + [A_{m}, B_{m}] = 0 \eqno (46a)
$$
$$
A_{mt} - C_{mx} + [A_{m}, C_{m}] = 0 \eqno (46b)
$$
$$
B_{mt} - C_{my} + [B_{m}, C_{m}] = 0. \eqno (46c)
$$
From (46a) we get
$$
k_{y} - m_{3x} + \sigma m_{1} - \tau m_{2} = 0 \eqno(47a)
$$
$$
\sigma_{y} - m_{2x} + \tau m_{3} - km_{1} =0   \eqno(47b)
$$
$$
\tau_{y} - m_{1x} + \beta (km_{2} - \sigma m_{3}) =0. \eqno(47c)
$$

The mM-LXII  equation (47a), we can rewrite in form
$$
k_{y} - m_{3x} = \beta {\bf e}_{3} \cdot ({\bf e}_{3x}\wedge {\bf e}_{3y})
 \eqno(48a)
$$
$$
\sigma_{y} - m_{2x} = \beta {\bf e}_{2}\cdot ({\bf e}_{2x}\wedge {\bf e}_{2y})
 \eqno(48b)
$$
$$
\tau_{y} - m_{1x}  = {\bf e}_{1} \cdot ({\bf e}_{1x}\wedge {\bf e}_{1y}).
\eqno(48c)
$$

Also from (14) we get
$$
k_{t} - \omega_{3x} = \sigma \omega_{1} - \tau \omega_{2} =
+\beta {\bf e}_{3}({\bf e}_{3x}\wedge {\bf e}_{3t})
\eqno(49a)
$$
$$
\tau_{t} - \omega_{1x} = \beta (k \omega_{2} - \sigma \omega_{3}) =
- {\bf e}_{1} \cdot ({\bf e}_{1x}\wedge {\bf e}_{1t})
 \eqno (49b)
$$
$$
\sigma_t- w_{2x}=\tau w_3-k w_1= \beta {\bf e}_{2} \cdot
({\bf e}_{2x} \wedge {\bf e}_{2t} \eqno(49c)
$$
and
$$
m_{1t} - \omega_{1y}= - \beta (m_{3} \omega_{2} - m_{2} \omega_{3}) =
e_1 \cdot (e_{1y} \wedge e_{1t})  \eqno (50a)
$$
$$
m_{2t} - \omega_{2y} =- m_{1} \omega_{3} + m_{3} \omega_{1} =
\beta e_2 \cdot (e_{2y}\wedge e_{2t})  \eqno (50b)
$$
$$
m_{3t} - \omega_{3y} =- m_{2} \omega_{1} + m_{1} \omega_{2} =
\beta e_3\cdot(e_{2y}\wedge e_{3t}). \eqno (50c)
$$

\subsection{On the topological invariants}

It is interesting to note that theM-LXII equations allows the following
integrals of motion
$$
K_{1} = \int \int (\kappa m_2 +\sigma m_3 dxdy, \quad
K_{2} = \int \int (-\tau m_3+k m_1) dxdy, \quad
K_{3} = \int \int (\tau m_2 -\sigma m_1 )dxdy \eqno(51a)
$$
or
$$
K_{1} = \int \int {\bf e}_1({\bf e}_{1x} \wedge {\bf e}_{1y})dxdy \eqno(51b)
$$
$$
K_{2} = \int \int {\bf e}_2({\bf e}_{2x} \wedge {\bf e}_{2y})dxdy
 \eqno(51c)
$$
$$
K_{3} = \int \int {\bf e}_3({\bf e}_{3x} \wedge {\bf e}_{3y})dxdy
 \eqno(51d)
$$

So we have the following three topological invariants
$$
Q_{1} =\frac{1}{4\pi} \int \int {\bf e}_{1}\cdot ({\bf e}_{1x}\wedge {\bf e}_{1y})dxdy
\eqno(52a)
$$
$$
Q_{2} =\frac{1}{4\pi} \int \int {\bf e}_{2}\cdot ({\bf e}_{2x}\wedge {\bf e}_{2y})dxdy
\eqno(52b)
$$
$$
Q_{3} =\frac{1}{4\pi} \int \int {\bf e}_{3}\cdot ({\bf e}_{3x}\wedge {\bf e}_{3y})dxdy
\eqno(52c)
$$

\subsection{The Lax   representation  of the mM-LXI equation}

To find the Lax  representation (LR) of the mM-LXI equation (44), we rewrite
it in the following matrix form
$$
\hat e_{1x}=
k \hat e_{2}   -\sigma \hat e_{3} \eqno(53a)
$$
$$
\hat e_{2x}=-\beta \hat e_{1} +\tau \hat e_{3} \eqno(53b)
$$
$$
\hat e_{3x}=\beta\sigma \hat e_{1}-\tau \hat e_{2}\eqno(53c)
$$
$$
\hat  e_{1y} =m_3 \hat e_{2} -m_2 \hat e_{3}
\eqno(54a)
$$
$$
\hat  e_{2y} =-\beta m_3 \hat e_{1} +m_1 \hat e_{3}
\eqno(54b)
$$
$$
\hat  e_{3y} =  \beta m_2\hat e_{1}-m_1 \hat e_{2}  \eqno(54c)
$$
$$
\hat e_{1t}=\omega_{3} \hat e_{2}-\omega_{2} \hat e_{3} \eqno(55a)
$$
$$
\hat e_{2t}=-\beta  \omega_{3} \hat e_{1} + \omega_1 \hat e_{3} \eqno(55b)
$$
$$
\hat e_{3t}=\beta\omega_{2} \hat e_{1}- \omega_{1} \hat e_{2}\eqno(55c)
$$
where $ \hat e_j $ are given by (14).
Equations (53)-(55) we can rewrite in the form (below we put $ \beta = 1 $)
$$
[\sigma_{3}, U] =k \sigma_{2}  -\sigma \sigma_{1} \eqno(56a)
$$
$$
[\sigma_{2},U] =-\beta k \sigma_{3}+\tau \sigma_{1} \eqno(56b)
$$
$$
[\sigma_{1},U] =\beta\sigma \sigma_{3}-\tau\sigma_{2}\eqno(56c)
$$
$$
[\sigma_{3},V]= m_{3} \sigma_{2} -m_{2} \sigma_{1}
\eqno(57a)
$$
$$
[\sigma_{2},V] =-\beta m_{3}  \sigma_{3} -m_{1} \sigma_{1} \eqno(57b)
$$
$$
[\sigma_{1},V] =\beta m_{2}\sigma_{3} - m_{1} \sigma_{2}  \eqno(57c)
$$
$$
[\sigma_{3}, W] = \omega_{3}\sigma_{2} - \omega_{2} \sigma_{1}\eqno(58a)
$$
$$
[\sigma_{2},W] = -\beta\omega_{3}\sigma_{3} +\omega_{1} \sigma_{1}\eqno(58b)
$$
$$
[\sigma_{1},W] =\beta \omega_{2} \sigma_{3} - \omega_{1} \sigma_{2}\eqno(58c)
$$
where
$$
U =  g_{x}g^{-1},\quad V =  g_{y}g^{-1}, \quad W =  g_{t}g^{-1}.  \eqno(59)
$$
Hence we get
$$
U =
\frac{1}{2i}
\left ( \begin{array}{cc}
\tau  & k-i\sigma \\
\beta(k+i\sigma) & -\tau
\end{array} \right)
\eqno(60a)
$$
$$
V=
\frac{1}{2i}
\left ( \begin{array}{cc}
m_1  & m_3 - im_2 \\
\beta(m_3 +im_2)  &  -m_1
\end{array} \right)
\eqno(60b)
$$
$$
W=
\frac{1}{2i}
\left ( \begin{array}{cc}
\omega_1  &  w_3 - iw_2 \\
\beta(w_3 + iw_1)  & -w_1
\end{array} \right).
\eqno(60c)
$$

Thus the matrix-function  $g$ satisfies the equations
$$
g_{x}=Ug,\quad g_{y}=Vg, \quad  g_{t}=Wg.   \eqno(61)
$$
This equation is the LR of the mM-LXI equation. Apropos as $\sigma = 0$
the equation (61) is the LR of the M-LXI equation (29). From the
compatibility condition of the equations (61) we get the new form of the
mM-LXII equation (46)
$$
U_{y}-V_{x}+[U,V]=0 \eqno(62a)
$$
$$
U_{t}-W_{x}+[U,W]=0 \eqno(62b)
$$
$$
V_{t}-W_{y}+[V,W]=0 \eqno(62c)
$$

\subsection{The M-III equation}

In this subsection,  we want establish the  connection between
the mM-LXI (44), mM-LXII (46) and M-III  (5)
equations. To this purpose, as above,  we identify $ {\bf  S}$
with ${\bf e}_{1}$, i.e. works the identity (39).
Then the identifying variables for the M-III equation (5) are given by
$$
m_{1}=u+\partial_{x}^{-1}\tau_{y} \eqno(63a)
$$
$$
m_{2}=\frac{1}{k}(u_x+\sigma m_3)  \eqno (63b)
$$
$$
m_{3} =\partial_x^{-1}(k_y+\sigma m_1-\tau m_2) \eqno(63c)
$$
and
$$
\omega_{1} = \frac{1}{k}[\sigma_{t}-\omega_{2x}+\tau\omega_{3}]
                 \eqno (64a)
$$
$$
\omega_{2}= -m_{3x}-\tau m_2+u \sigma+2l(cl+d)m_2- 4cv\sigma
\eqno (64b)
$$
$$
\omega_{3}= m_{2x}-\tau m_3-uk+2l(cl+d)m_3-4cvk
\eqno (64c)
$$

\subsection{The L-equivalent of the M-III equation}

Return to the  function $q$. Let   this function has the form
$$
q =\frac{k^2+\sigma^2}{2(2cl+d)}exp i\{2l(cl+d)x-\partial^{-1}_x \tau\}
\eqno(65)
$$

Then $q$ again satisfies is the equation (23). So, we have again shown that the M-III
equation  (5)
and   the equation (23) are L-equivalent to each other. It is remarkable  that
this result is consistent with the other result namely that
these equations are G-equivalent (gauge equivalent) to each other [3].

\section{Geometry of curves and Bilinear representation of the M-III equation}

In this section we  establish self-coordination
of the our geometrical formalism that presented above with the other
powerful tool of soliton theory - the Hirota's bilinear method.
We demonstrate our idea in example the M-III equation (5). For the curve
we take the mM-LXI and mM-LXII equations.  Usually,
for the spin vector ${\bf S}=(S_{1},S_{2},S_{3})$ takes the following transformation
$$
S^{+} = S_{1}+iS_{2}= \frac{2\bar f g}{\Lambda}, \quad S_{3} = \frac{\bar f f
- \bar g g}{\Lambda}, \quad \Lambda=\bar f f+\bar g g.  \eqno(66)
$$
Also in this section, we assume the (39) is holds.
In [4] was shown that for the mM-LXI equation (44) is correct the following
representation
$$
e^{+}_{1} = \frac{2\bar f g}{\Lambda}, \quad e_{13} = \frac{\bar f f
- \bar g g}{\Lambda} \eqno(67a)
$$
$$
e^{+}_{2} = i\frac{\bar f^{2} +\bar g^{2}}{\Lambda}, \quad
e_{23} = i\frac{fg- \bar f \bar g}{\Lambda} \eqno(67b)
$$
$$
e^{+}_{3} = \frac{\bar f^{2} - \bar g^{2}}{\Lambda}, \quad
e_{33} = - \frac{fg + \bar f \bar g}{\Lambda}  \eqno(67c)
$$
and
$$
k =-i\frac{D_{x}(g \circ f - \bar g\circ \bar f )}{\Lambda},
\quad \sigma = -\frac{D_{x}(g \circ f + \bar g\circ \bar f )}{\Lambda},
\quad
\tau=-i\frac{D_{x}(\bar f \circ f + \bar g\circ g )}{\Lambda} \eqno(68a)
$$
$$
m_{1}=-i\frac{D_{y}(\bar f \circ f + \bar g\circ g )}{\Lambda}
\quad
m_{2}=-\frac{D_{y}(g \circ f + \bar g\circ \bar f)}{\Lambda}
\quad m_{3}=-i\frac{D_{y}(g\circ f - \bar g\circ \bar f)}{\Lambda}
  \eqno(68b)
$$

Here ${\bf e}_{j} = (e_{j1}, e_{j2}, e_{j3}), \quad e_{j}^{\pm} = e_{j1}
\pm ie_{j2}$. Now we take
$$
\tau = 0 \quad m_{1} =  u.  \eqno(69)
$$
Then hence and from (68) we have
$$
D_{x}(\bar f \circ f + \bar g\circ g ) = 0 \eqno(70a)
$$
$$
u = - i\frac{D_{y}(\bar f \circ f + \bar g\circ g )}{\Lambda} \eqno(70b)
$$
So for the M-III equation (5) we have the bilinear representation (66) and (70).
This representation allows construct the bilinear form of the M-III equation
which is left in future.

\section{The M-III equation as the integrable particular case
  of the M-0 equation}

Consider the (2+1)-dimensional M-0 equation [1]
$$
{\bf S}_{t} = a_{12} {\bf e}_{2} + a_{13}{\bf e}_{3}, \quad
{\bf S}_{x} = b_{12} {\bf e}_{2} + b_{13}{\bf e}_{3}, \quad
{\bf S}_{y} = c_{12} {\bf e}_{2} + c_{13}{\bf e}_{3} \eqno(71)
$$
where
$$
{\bf e}_{2} = \frac{c_{13}}{\triangle}{\bf S}_{x} -
\frac{b_{13}}{\triangle}{\bf S}_{y}, \quad
{\bf e}_{3} = -\frac{c_{12}}{\triangle}{\bf S}_{x} +
\frac{b_{12}}{\triangle}{\bf S}_{y}, \quad \triangle =
b_{12}c_{13}- b_{13}c_{12}. \eqno(72)
$$
All known spin systems (integrable and nonintegrable) in 2+1 dimensions
are the particular reductions of the M-0 equation (71). In particular,
the M-III equation (5) is the integrable reduction of equation (71). In this case,
we have
$$
a_{12} = \omega_{3}, \quad a_{13}= -\omega_{2}, \quad b_{12}= k, \quad
b_{13}=  -\sigma, \quad
 c_{12}= m_{3}, \quad c_{13}= -m_{2}.  \eqno(73)
$$
Sometimes we use the following form of the M-0 equation [1]
$$
{\bf S}_{t} = d_{2} {\bf S}_{x} + d_{3}{\bf S}_{y} \eqno(74)
$$
with
$$
d_{2}= \frac{a_{12}c_{13}-a_{13}c_{12}}{\triangle}, \quad
d_{3}= \frac{a_{12}b_{13}-a_{13}b_{12}}{\triangle}. \eqno(75)
$$

\section{The  equation for $\lambda$}

Let us consider the M-LXVII equation in he form  (12)
for the  case $q=p=0$.  We have
$$
\left ( \begin{array}{c}
{\bf e}_{1} \\
{\bf e}_{2} \\
{\bf e}_{3}
\end{array} \right)_{x} =
\left ( \begin{array}{ccc}
0      &   0                      & 0 \\
0      & 0                        & -2(c\lambda^{2}+d\lambda)  \\
0      & 2(c\lambda^{2}+d\lambda) & 0
\end{array} \right)
\left ( \begin{array}{ccc}
{\bf e}_{1} \\
{\bf e}_{2} \\
{\bf e}_{3}
\end{array} \right) \eqno(76a)
$$
$$
\left ( \begin{array}{ccc}
{\bf e}_{1} \\
{\bf e}_{2} \\
{\bf e}_{3}
\end{array} \right)_{t}= 2(c\lambda^{2} + d \lambda)
\left ( \begin{array}{ccc}
{\bf e}_{1} \\
{\bf e}_{2} \\
{\bf e}_{3}
\end{array} \right)_{y}  \eqno(76b)
$$
Hence we get
$$
\lambda_{t} =  2(c\lambda^{2} + d \lambda)\lambda_{y}  \eqno(77)
$$
So for the M-III equation (5) the spectral parameter satisfies the equation
(77), i.e. in this case we have a nonisospectral problem. From (77) we obtain \\
1) for the M-I equation (2)
$$
\lambda_{t} =  2 \lambda \lambda_{y}.  \eqno(78)
$$
2) for the M-II equation (3)
$$
\lambda_{t} =  2c\lambda^{2}\lambda_{y}.  \eqno(79)
$$
Now consider the general form of such equations
$$
\lambda_{t} =  k \lambda^{n}\lambda_{y}, \quad k=const.  \eqno(80)
$$
This equation has the following solution
$$
\lambda   = ( \frac{y+c}{a-kt})^{\frac{1}{n}}  \eqno(81)
$$
where a(c) is real (complex) constant.

\section{Conclusion}

In this paper, we have used a geometrical approach pioneered by Lakshmanan
 to analyze the connection between differential geometry of curves and
spin systems to establish such connection with the integrable (2+1)-
dimensional spin system- called M-III equation. Simultaneously our approach
permit construct the corresponding L-equivalent of the given spin system and
Lax representation of it. Some other consequences of this geometrical formalism
are  presented.

 {99}
\end{document}